\makeindex 
\documentclass[leqno]{conm-p-l}
\usepackage{epsfig,amssymb} 
\newcommand{\R}{\mathbb{R}}
\newcommand{\C}{\mathbb{C}}
\newcommand{\cM}{\mathcal{M}}

\renewcommand{\qed}{$\blacksquare$}
\newcommand{\cC}{\mathcal{C}}
\newtheorem{ex}{Example}

\newtheorem{thm}{Theorem}[ex]

\newenvironment{coxitem}
{\begin{list}
{{$\bullet$}}{\setlength{\leftmargin}{18pt}\setlength{\topsep}{0pt}}}
{\end{list}}

\newcounter{cxlist}

\newenvironment{coxenum}
{\begin{list}
{(\alph{cxlist})}{\usecounter{cxlist}\setlength{\topsep}{0pt}
\setlength{\leftmargin}{18pt}}}
{\end{list}}

\begin{document}

\title[The Degree of an equation obtained by elimination]{On the
    determination of the degree of an equation obtained by
    elimination}

\author{Ferdinand Minding}

\begin{abstract} In 1841, Ferdinand Minding published \emph{Ueber die
Bestimmung des Grades einer durch Elimination hervorgehenden
Gleichung} in Crelle's Journal.  His main theorem
represents the first (implicit) appearance of the BKK bound (and even
some of its extensions) in the mathematical literature.  We present
an English translation of this paper and a commentary which explains
how Minding's formula relates to the mixed area of lattice polygons.
\end{abstract}


\vbox{

\maketitle

\vglue-.15in

\centerline{Translation and Commentary by D.\ Cox and J.\ M.\ Rojas}

\smallskip

\setcounter{page}{351}

\section{Translation of Minding's Paper \cite{minding}}
{}From two given algebraic equations involving two unknowns, there are
indeed different methods for deriving a new equation containing only
one unknown.  Frequently one wishes to know only the degree of this
final equation, without writing down the equation itself.  If I am not
mistaken, there has until now been no rule given for finding this
degree with the ease appropriate to so elementary an object.  One
knows, to be sure, that if the equations are respectively of degrees
$h$ and $k$, i.e., if the highest sum of the exponents of terms
amounts to $h$ in one and $k$ in the other, then the desired degree
can equal at most the product $hk$.  This however is only a bound,
from which the actual degree often deviates.  In order to find the
actual degree, it is therefore necessary to assess the true form of the
final equation.\index{support}  We have the following equations:
\begin{align}
{f}(x,y)&=A_0y^m + A_1y^{m-1} + A_2y^{m-2}+\dots+A_{m-1}y+A_m = 0, \\
\theta(x,y)&=B_0y^n+B_1y^{n-1}+B_
2y^{n-2}+\dots+B_{n-1}y+B_{n} = 0,
\end{align}
in which the letters $A$ and $B$ with subscripts stand for arbitrary
polynomials in $x$.  We solve the equation (2) for $y$, denote its roots
by ${y_1, y_2,\dots, y_n,}$ and form the product
\begin {equation}                        
P = f (x, y_1) \cdot f (x, y_2) \cdots f(x,y_n).
\end{equation}
Then ${B_0^m}\cdot{P}$ is a polynomial in $x$, and if one sets
\begin{equation}
{B_0^m}\cdot{P}=\psi(x),
\end{equation}
then
\begin{equation}
\psi(x)=0
\end{equation}} 
\newpage
\noindent
is the desired final equation.

In order to show that $B_0^m\cdot{P}$ is a polynomial, first note that
$P$ is a rational function of $x$, in which, if the functions of $y_1,
y_2,\dots, y_n$ occurring symmetrically therein are expressed in terms
of $x$ using equation (2), the denominator can only contain a
power of $B_0$.  Designate any one of these symmetric functions by
\begin{equation*}
S = {y_1^{m_1}} {y_2^{m_2}}{y_3^{m_3}}\cdots{y_n^{m_n}}+\cdots,
\end{equation*}
where the following terms indicated by the dots are obtained from the
first through mixing $y_1, y_2,\dots, y_n$, and none of the exponents
$m_1, m_2, \dots, m_n$ can be greater than $m$.  Now set:         
\begin{equation*}
{S_1}=\frac{1}{y_1^{m-{m_1}}\cdot{y}_2^{m-{m_2}}
\cdots{y}_n^{m-{m_n}}}+\cdots,
\end{equation*}
therefore
\begin{equation*}
{S} = 
(y_1 y_2\cdots y_n)^m {S_1}=(-1)^{mn}\frac{B_n^m}{B_0^m}{S_1}.
\end{equation*}
The symbol $S_1$ denotes a symmetric polynomial of the reciprocal roots
of (2).  Therefore the value of $S_1$ is a rational function, whose
denominator can have only a power of $B_n$; so we set
$S_1=\frac{Z}{B_{n}^\lambda},$ where $Z$ is a polynomial in $x$ or
more accurately a polynomial function of the polynomials $B_0, B_1,
\dots, {B}_n$, and $\lambda$ is a positive integer.  It follows that
\begin{equation*}
S=(-1)^{mn}\frac{B_{n}^{m}\cdot{Z}}{B_{0}^{m}\cdot{B}_{n}^{\lambda}}.
\end{equation*}
Since $S$ obviously can only have a power of $B_0$ in its denominator,
$B_{n}^{\lambda}$ must appear in the above numerator; therefore
$B_{0}^{m}{S}$ and consequently $B_{0}^{m}{P}=\psi(x)$ are
polynomials.

When equation (1) is solved for $y$, we denote the roots by $\eta_{1},
\eta_{2},\dots,\eta_{m}$.  Set
\begin{equation*}
{Q}=\theta(x,\eta_{1})\cdot\theta(x,\eta_{2})\cdots\theta(x,\eta_{m})
\end{equation*}
and notice that
\begin{equation*}
\theta(x,\eta_1)=B_0(\eta_1-y_1)(\eta_1-y_2)\cdots(\eta_1-y_n),\ 
\text{etc.},
\end{equation*}
so that
\begin{multline*}
Q=B_0(\eta_1-y_1)\cdots
(\eta_1-y_n)\times{B}_0(\eta_2-y_1)\cdots(\eta_2-y_n)\times\cdots\\
\cdots\times{B}_0(\eta_m-y_1)\cdots(\eta_m-y_n),
\end{multline*}
and because
\begin{equation*}
A_{0}(y_1-\eta_1)(y_1-\eta_2)\cdots(y_1-\eta_m)=f(x,y_1),\ \text{etc.},
\end{equation*}
it follows that:
\begin{equation}
A_{0}^{n}Q=(-1)^{mn}B_{0}^{m}P.
\end{equation} 

This proves that $\psi(x)=0$ is the required final equation.  Namely
for every value assumed in this problem by $x$, $P = 0$ (as well as $Q
= 0$) follows necessarily, thus $\psi(x)=0$.  Furthermore, should this
equation contain a superfluous factor, then neither $P = 0$
nor $Q = 0$ for such a $\psi(x) = 0$; but then, because of (4) and
(6), $A_0$ and $B_0$ would vanish together, which in general\footnote{ 
{\bf Translator's Note:} In 
Minding' paper, ``in general'' does \emph{not} imply complete generality. 
In modern parlance, ``in general'' should be interpreted as  
``generically''.\index{genericity} } is not possible.  If $A_0$ and 
$B_0$ have a common factor in a
particular case, then the polynomial $\psi(x)$ is always divisible by
this factor, because, as is easily seen, it always is of the form
$\psi(x)=A_0U+{B_0}V$, in which $U$ and $V$ are polynomials.  Because
one can always remove such a case through an infinitely small change
of the coefficients,\index{genericity} and to be sure without changing the 
degree, so it
follows, that in no case\footnote{{\bf Translator's Note:} Here some
care is required.  One in general needs to perturb more coefficients
than just those of $A_0$ and $B_0$. For example, taking
$f=ay^4+x^2y^3+x^3y^2+b$ and $\theta =
cy^4+x^2y^3+x^2y^2+d$, even with $a,b,c,d$ generic, yields a system
which always has one root at infinity,\index{root at infinity} thus 
violating Minding's formula.  See also Section \ref{sub:gen}. }
does the equation $\psi(x)=0$ have a root
which doesn't come from a solution of our problem.

The degree of the polynomial $\psi(x)$ reveals itself in the following
way.  We have
\begin{equation*}                 
\psi(x)=B_{0}^{m}\cdot{f}(x,y_1)\cdot{f}(x,y_2)\dots{f}(x,y_n).
\end{equation*}
One develops the roots ${y}_1, {y}_2,\dots,{y}_n$ of equation (2)
according to decreasing powers of $x$\index{Puiseux series} and sets 
the series in place of each root in the above expression.\footnote{ 
{\bf Translator's Note:} By ``decreasing powers'', Minding means
a series of the form 
\[
y(x) = cx^{h} + dx^{h-\alpha_{1}} + ex^{h-\alpha_{2}} + \cdots, \quad
0 < \alpha_{1} < \alpha_{2} < \cdots,\quad h,\alpha_{1}, \alpha_{2},
\ldots \in \mathbb{Q}.
\]
Minding calls $cx^{h}$ the \emph{first term} and $h$ the \emph{highest
exponent}. }
Then all fractional and negative powers of $x$ will mutually cancel
themselves, and we obtain the polynomial $\psi(x)$ unchanged.  Since
only the degree of $\psi(x)$ is desired, one instead replaces each
series with only its first term, which for ${y}_1, {y}_2, \dots,{y}_n$
are $c_1{x}^{h_1}, {c}_2{x}^{h_2},\dots, {c}_n{x}^{h_n}$.  The
procedure, through which the series and in particular the highest
exponents $h_1, {h}_2, \dots, {h}_n$ or the degree of the roots is
found, is sufficiently known; one can compare for example
\textit{Lacroix Traite} page 223 of the first edition, where the
development according to increasing powers is shown.  One determines
from this the highest exponent of $x$ in each of the functions $f(x,
{c}_1{x}^{h_1}),\ f(x, {c}_2{x}^{h_2}), \dots $ or the \emph{degree}
of the functions $f(x, y_1),\ f(x, y_2), \dots$, which will be denoted
by ${k}_1, {k}_2,\dots, {k}_n$.  These can be integers or fractions,
but can never be negative since $A_m$ has degree at least $0$.\footnote{ 
{\bf Translator's Note:} One needs to forbid negative exponents 
in \emph{all} $A_i$, not just $A_m$.}  
Finally, we denote the degree of $B_0$ by $b$, so that
\begin{equation}
mb+{k}_1 +  {k}_2 +{k}_ 3 + \dots + {k}_n 
\end{equation}
is necessarily an integer, which gives the highest exponent of
$\psi(x)$ or the sought-for degree of the final equation.\index{Bernstein's 
Theorem} In special
cases, one can consider the values of $c_1, {c}_2,\dots, {c}_n$ in
order to see, whether the coefficient of the highest term in one of
the factors $f(x,y_1), \dots$ of $\psi(x)$, consequently of $\psi(x)$
itself, perhaps becomes identically zero, and in such a case, one will
be compelled to include the subsequent terms of the series in the
calculation; however it is not necessary to continue this suggestion
any further.  Rather it is clear, that in general the
above-mentioned value (7) represents the actual degree of the
polynomial $\psi(x)$.\index{mixed area} 

Let the following two equations be given,\index{support} in which the symbol
$(x^{\mu})$ denotes a polynomial of degree $\mu$ in $x$:
\begin{align*}
f(x,y)& = (x^2)y^4+(x^2)y^3+(x^4)y^2+(x^5)y+(x^5) = 0,\\
\theta(x, y)& = (x^8)y^5+(x^6)y^4+(x^9)y^3+(x^4)y^2+(x^3)y+(x^4) = 0.
\end{align*}
These equations are of the 6th and 13th degree, so that the degree of
the final equation is not higher than $6\cdot13 = 78$. To find it
exactly, one calculates the degree of the roots $y$ of $\theta(x, y) =
0$; one finds at once that $h_1 = h_2 = \frac{1}{2}, h_3 = h_4 = h_5 =
{-} \frac{5}{3}$.  It follows that the degrees of $f(x, y_1), \dots$ are 
namely $k_1 = k_2 = \frac{11}{2}$ and $k_3 = k_4 = k_5 = 5$; furthermore
$B_0 = (x^8)$, so that $b = 8$, and $m = 4$.  Thus the degree of the
final equation is $mb + k_1 + k_2 + k_ 3 + k_4 +k_5 = 4\cdot 8 + 11 +
15 =58$.

If one writes the given equations in terms of $x$ instead of $y$, in
order to find the degree in $y$ of the final equation using the 
above-mentioned rule, one does not always find the same value for this
as for the previous degree.  To explain this circumstance, one must
note, that the final equation in $x$ only gives the finite values of
$x$, which suffice to produce zero in both equations.  If the final
equation in $y$ has a higher degree than that in $x$, then it is
necessary that some root of the equation in $y$ belong to an infinite
value of $x$.\index{root at infinity} It is always easy to make this value 
finite through an infinitely small change of the 
coefficients\index{genericity} of one of the 
existing equations and to remove the inequality of the degrees of the final
equations.  Namely, let the equations (1) and (2) be written as
follows:
\begin{align*}
f(x,y)& = \alpha_0{x}^\mu+\alpha_1{x}^{\mu-1}+ \dots+\alpha_\mu = 0,\\     
\theta(x, y)& = \beta_0{x}^\nu+\beta_1{x}^{\nu-1}+ \dots+\beta_\nu = 0,
\end{align*}
where $\alpha_0, \alpha_1, \dots, \beta_0, \dots,\beta_\nu$ are
polynomials in $y$.  If now neither $A_0$ has a common factor with
$B_0$, nor $\alpha_0$ with $\beta_0$, then neither can $y$ have an
infinite value for finite $x$, nor infinite $x$ for finite
$y$;\footnote{ {\bf Translator's Note:} From a modern point of view,
the roots lie on a toric surface such that roots will not have finite
$x$ or $y$ coordinates when they lie on curves corresponding to edges
that are neither vertical nor horizontal in the Newton polygon.  It is
amazing that Minding knew this.\index{root at infinity}} 
consequently there can then be no difference between the degrees of
the final equations in $x$ and in $y$.  If a joint factor between
$A_0$ and $B_0$ or $\alpha_0$ and $\beta_0$ is present, one needs
therefore only to change a coefficient in $A_0$ or one in $\alpha_0$
in order for the final equations in $x$ and $y$ to have the same
degree.  If immediately afterwards one sets these changes to zero, then
one can examine the coefficients of the highest term of the equation in
order to decide how many values of $x$ and how many of $y$ become
infinite, and finally how many finite solutions of the equations are
present.  This version of the calculation is however unnecessary if one
uses the stated rule properly.  Consider for example the following
equations:
\begin{align*}
f(x, y) & = (a+bx^2)y^4+(c+ex)y^2+gx^3y+h +kx^2+lx^3 = 0, \\
\theta(x,y) & =  \beta{x}^5y^2+(\gamma+\delta{x}^2)y+ 
\lambda+\mu{x}^4 = 0, 
\end{align*}
or ordered with respect to $x$:
\begin{align*}
f(x,y) & = (l+gy)x^3 + (k+by^4)x^2+ey^2x +h + cy^2 +ay^4 = 0, \\
\theta(x, y) & = \beta{y}^2x^5+ \mu{x}^4 +\delta{yx}^2 +\lambda +
\gamma{y} = 0.
\end{align*}
Here $A_0 = a+ bx^2, B_0 =\beta{x}^5, \alpha_0 = l+gy, \beta_0 =
\beta{y}^2$; consequently, unless $a = 0$ or ${l} = 0$, $A_0$ and
$B_0$ as well as $\alpha_0$ and $\beta_0$ have no common factor, so
that the degrees of the final equation in $x$ and $y$ are both found
to be $=26$.  However, if one sets $a = 0$ and $l = 0$ at the same
time, and calculates the degree of the final equations, then one finds
$25$ for the equation in $x$, and $24$ for that in $y$.  Through the
vanishing of $a$ and $l$, two roots of the previous final equation in
$y$ and one of the previous final equation in $x$ become therefore
infinite; at the same time however also the new final equation in $x$
is divible by ${x}^2$, the common factor of $A_0$ and $B_0$, so
likewise the new final equation in $y$ is divisible by $y$, the common
factor of $\alpha_0 $ and $\beta_0$.  Of the $26$ finite solutions,
which come from the initial equations, $23$ in general remain finite
when $a$ and $l$ vanish; the three remaining become $x_{24} = 0$,
$y_{24} = \infty$; $x_{25} = 0$, $y_{25} = \infty$; $x_{26} = \infty$,
$y_{26} = 0$.  One sees, how here the desired number of finite
solutions is found through repeated use of the stated rule and
comparison of the results.

\bigskip

{\bf Postscript.} After completing this project, a new work came to my
attention: \emph{System der Algebra von Dr.\ P.J.E.\ Finck, Professor
zu Stra{\ss}burg; Leipzig bei Barth, 1841}, which on page 405
indicates a much more concise rule for the degree of the final
equation than the one mentioned in this project.  The rule (whose 
proof is given below) is the following:\ if all coefficients $A_0, A_1,
\dots, A_m$ of equation (1) have degree $m^\prime$ and all
coefficients $B_0, B_1, \dots, B_n$ of equation (2) have degree
$n^\prime$, then the degree of the polynomial $\psi(x)$, which gives
the final equation, is the following:\ ${m}n^\prime + n m^\prime$.  The
proof of this theorem, which fills two pages of the text, follows very
easily from the above, for in this case all roots $y$ of equation (2)
have degree $0$, that is, $h_1 = h_2 = \dots = h_n = 0$, so that $k_1
= k_2 = \dots = k_n = m^\prime$; at the same time, $b = n^\prime$
because $B_0$ has degree $n^\prime$.  Consequently the degree of the
final equation is $m b + k_1 + k_2 + \dots + k_n =m n^\prime + n
m^\prime$, qed.  If this rule is used in cases in which the
coefficients are of unequal degree, then the result is no longer
dependable because it assumes equality of every degree.  This can
result in extraneous roots.  On the other hand, in the rule stated in
the present treatise, in order to find the degree of $\psi(x)$, the
degrees of the coefficients are used correctly in the calculation no matter 
how they are given. 

\section{Commentary on Minding's Paper}

The main result of Minding's paper \cite{minding} is the formula

\smallskip\vskip2pt

\noindent (7) \hfill $mb+{k}_1 +  {k}_2 +{k}_ 3 + \dots +
	  {k}_n$\hfill\mbox{}  

\smallskip\vskip2pt

\noindent for the degree of the \emph{final equation}
(\emph{Endgleichung}) obtained by eliminating $y$ from equations (1)
and (2).  Minding notes that this degree is the number of solutions of
the equations, at least when the coefficients are generic.  On the 
other hand, the BKK bound states that if $P_1$ and $P_2$ are the Newton
polygons of the polynomials in (1) and (2), then for generic
coefficients, the number of solutions with all coordinates nonzero is
the mixed area\index{mixed area} 
$\mathcal{M}(P_1,P_2)$ \cite{bernie,bz,rojas}.  The
goal of this commentary to explain why these formulas give the same
answer for the types of polynomials considered by Minding, and how
Minding's result also foreshadowed later extensions of the BKK bound
to counting all complex solutions.  We will also clarify the implicit
genericity assumptions.

Although (7) seems far removed from mixed areas, recall that Minding
expresses the roots $y$ of $\theta(x,y)$ as certain fractional power
series, now called \emph{Puiseux series}.\index{Puiseux series}  
The technique for computing
these series is due to Newton and is where he introduced the Newton
polygon.  Once we understand how this works, we will see that (7)
reduces to standard formulas for the mixed area of lattice polygons.

An introduction to Puiseux series appears in \cite[Ch.\ IV]{walker}
\index{Puiseux series} 
and the properties we need about mixed areas can be found in
\cite[Ch.\ 7]{clo} and \cite[Sec.\ 7--8]{rojas}.\index{mixed area}  
For consistency, we will use the notation of \cite{rojas} for mixed areas.

\subsection{Puiseux Series and Newton Polygons} 
\label{sub:puiseux} 
If $y_1.\dots,y_n$
are the roots of $\theta(x,y) = 0$, then standard resultant formulas
\cite[(1.4) of Ch.\ 3]{clo} show that the right-hand side of
\[
\psi(x) = B_0^m f(x,y_1) f(x,y_2) \cdots f(x,y_n)
\]
is the univariate resultant $\mathrm{Res}(\theta,f,y)$.  This gives
another proof that $\psi$ is the polynomial obtained by eliminating
$y$ from (1) and (2).

Minding's innovation was to express the $y_i$ in terms of $x$ using
Puiseux series
\begin{equation}
\label{puiseaux}
y(x) = cx^{h} + dx^{h-\alpha_{1}} + ex^{h-\alpha_{2}} + \cdots, \quad
0 < \alpha_{1} < \alpha_{2} < \cdots,
\end{equation}
where $h,\alpha_1,\alpha_2,\ldots \in \mathbb{Q}$.  To see how this
relates to Newton polygons, consider Minding's bivariate polynomial
\begin{equation}
\label{minth}
\theta(x, y) = (x^8)y^5+(x^6)y^4+(x^9)y^3+(x^4)y^2+(x^3)y+(x^4) = 0,
\end{equation}
where, following Minding, $(x^\mu)$ indicates a polynomial in $x$ of
degree $\mu$.  If all the coefficient polynomials $(x^\mu)$ are generic, 
then the Newton polygon\index{Newton polygon} $P_2$ of $\theta(x,y)$ is

\vglue.5in

\begin{equation}
\label{minpoly}
\end{equation}

\vglue-1.2in

\[
\epsfig{file=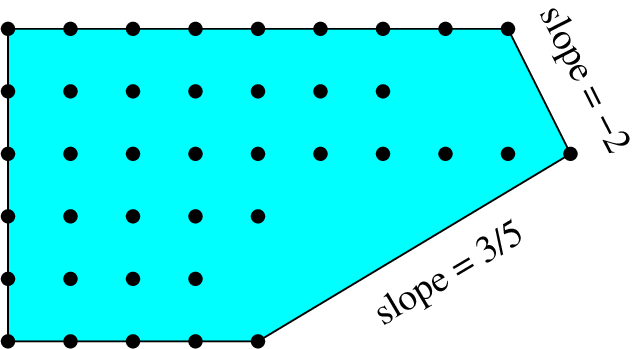,height=1.5in}
\]
To see how \eqref{puiseaux} relates to $P_2$, 
we write
\[
\theta(x,y) = \cdots + B_{8,5} x^8 y^5 + \dots + B_{9,3} x^9 y^3 +
\dots + B_{4,0} x^4.
\] 
If we substitute \eqref{puiseaux} into this and set the resulting
expression to zero, we obtain
\[
0 = \cdots + B_{8,5} x^8 (cx^{h} + \cdots)^5 + \dots + B_{9,3} x^9
(cx^{h} + \cdots)^3 + \dots + B_{4,0} x^4.
\]
The terms indicated by $\cdots$ in the parentheses have smaller degree
in $x$, so that the only way the highest degree $x$ terms can cancel
is for
\[
B_{8,5}c^5 x^{8+5h} + B_{9,3} c^3 x^{9+3h} = 0\quad \text{or}\quad
B_{9,3} c^3 x^{9+3h} + B_{4,0} x^4 = 0.
\]
The first equality implies that $h = \frac12$, and the (nonzero)
values of $c$ giving $B_{8,5}c^5 + B_{9,3} c^3 = 0$ yield two
solutions $y_1,y_2$ of $\theta(x,y) = 0$. The second equality above
implies that $h = -\frac53$, and the (nonzero) values of $c$ giving
$B_{9,3} c^3 + B_{4,0}= 0$ yield the remaining three solutions
$y_3,y_4,y_5$.  The exponents $\frac12$ and $-\frac53$ are the ones
mentioned by Minding in his paper.  Note that they are the negative
reciprocals of the slopes indicated in \eqref{minpoly}.

In general, assume that the Newton polygon of $\theta(x,y)$ has the
following form:

\vglue.4in

\begin{equation}
\label{genpoly}
\end{equation}

\vglue-1.15in

\[
\begin{picture}(210,120)
\put(0,5){\line(1,0){210}}
\put(5,0){\line(0,1){110}}
\thicklines
\put(5,5){\line(1,0){120}}
\put(125,5){\vector(2,1){40}}
\put(185,85){\vector(-1,1){20}}
\put(165,25){\vector(1,2){20}}
\put(5,105){\line(1,0){160}}
\put(5,5){\line(0,1){100}}
\put(180,95){$\gets\ v_1 = (n_1,m_1)$}
\put(150,10){$\gets\ v_r = (n_r,m_r)$}
\put(180,40){$\gets\ v_{r-1} = (n_{r-1},m_{r-1})$}
\put(185,70){$\vdots$}
\end{picture} 
\]
Note that the edges on the right are considered as vectors.
Then one can show that $\theta(x,y) = 0$ has $n$ solutions, where
\begin{equation}
\label{gensol}
\begin{aligned}
m_1\ \text{solutions have the form}\ y(x) &= c_{1,j}
x^{-\frac{n_1}{m_1}} + \cdots,\ j = 1,\dots,m_1,\\
m_2\ \text{solutions have the form}\ y(x) &= c_{2,j}
x^{-\frac{n_2}{m_2}} + \cdots,\ j = 1,\dots,m_2,\\
&\ \,\vdots \\
m_r\ \text{solutions have the form}\ y(x) &= c_{r,j}
x^{-\frac{n_r}{m_r}} + \cdots,\ j = 1,\dots,m_r. 
\end{aligned}
\end{equation}
and the dots $\cdots$ indicate lower order terms.\index{Puiseux series}  
The solutions involve the edges on the \emph{right-hand} side of the Newton 
polygon because Minding writes his Puiseux series with \emph{decreasing}
exponents. Note also that this $19^{\text{\underline{th}}}$ century
construction is a foreshadowing of general results for systems of
multivariate polynomial equations that form the core of what is now
called non-Archimedean amoeba theory or tropical algebraic geometry.

\subsection{Mixed Area of Lattice Polygons} Consider the standard
lattice $\mathbb{Z}^{2} \subset \mathbb{R}^{2}$.  The area of a
lattice polygon $P$ will be denoted $\mathrm{Area}(P)$, normalized so
that the unit square has area $1$.  Given two lattice polygons $P_{1},
P_{2}$, there are four methods to compute the mixed area  
$\mathcal{M}(P_{1},P_{2})$:\index{mixed area} 

\smallskip

\begin{coxitem}
\item (Definition) $\mathcal{M}(P_{1},P_{2})$ is the
  coefficient of $\lambda \mu$ in the expansion of the quantity 
  $\mathrm{Area}(\lambda P_{1} + \mu P_{2})$ as a polynomial in
  $\lambda,\mu$ \cite[Def.\ 7.0.25]{rojas}.
\item (Inclusion-Exclusion) $\mathcal{M}(P_{1},P_{2})\! =\!
  \mathrm{Area}(P_{1}+P_{2}) - \mathrm{Area}(P_{1}) -
  \mathrm{Area}(P_{2})$ \cite[Lem.\ 7.0.29]{rojas}.
\item (Mixed Subdivision) $\mathcal{M}(P_{1},P_{2})$ is the sum of the
  areas of the mixed cells in a mixed subdivision of $P_{1}+P_{2}$
  \cite[Lem.\ 7.0.29]{rojas}.
\item (Recursion) $\mathcal{M}(P_{1},P_{2}) = -\sum_{\mathcal{F}}
  \min_{u \in P_{1}} (u \cdot \nu_{\mathcal{F}})
  \mathrm{Length}(\mathcal{F})$.  The sum is over all edges
  $\mathcal{F}$ of $P_{2}$ where $\nu_{\mathcal{F}}$ is the
  primitive\footnote{Thus $\nu_{\mathcal{F}}$ is the 
first lattice point on the open ray generated by 
any inner normal of $\mathcal{F}$.} 
inward-pointing normal vector to $\mathcal{F}$ and
  $\mathrm{Length}(\mathcal{F})$ is the normalized length of
  $\mathcal{F}$ \cite[Thm.\ 4.12 of Ch.\ 7]{clo}.
\end{coxitem}

\smallskip

\noindent We will show that Minding's formula is equivalent to the
last two of these methods.  To make the relation to (7) easier to see, we
will write the recursion formula as
\begin{equation}
\label{mfor}
\mathcal{M}(P_{1},P_{2}) = \sum_{\mathcal{F}} \max_{u \in P_{1}}\, (u \cdot
\tilde{\nu}_{\mathcal{F}}).
\end{equation}
where as above, $\mathcal{F}$ is an edge of $P_{2}$ and
$\tilde{\nu}_{\mathcal{F}} = - \mathrm{Length}(\mathcal{F})\, 
\nu_{\mathcal{F}}$.

In practice, the vector $\tilde{\nu}_{\mathcal{F}}$ is easy to read
off from the polygon $P_{2}$.  Write $\mathcal{F}$ as a vector $v =
(n,m)$ with the polygon to the left of $v$.  This is the convention
used in \eqref{genpoly}.  Then one can easily show that
\begin{equation}
\label{tinu}
\tilde{\nu}_{\mathcal{F}} = (m,-n).
\end{equation}
So $\tilde{\nu}_{\mathcal{F}}$ is the outward-pointing normal vector
of $\mathcal{F}$ whose length equals that of $\mathcal{F}$.

\subsection{Minding's Formula and Recursion} 
\label{recursion}
We can now prove that Minding's 
formula (7) computes the mixed area of $P_1\!=$``Newton polygon of $f(x,y)$'' 
and $P_2\!=$``Newton polygon of $\theta(x,y)$'', assuming that both
polygons have the form \eqref{genpoly}. (See Section 
\ref{sub:mix}, just after Example 1, for a rigorous statement 
of the required conditions for $P_1$ and $P_2$.) 

In $P_2$, \eqref{tinu} shows the edges on the right in \eqref{genpoly}
give the following outward-pointing normal vectors: 
$\tilde{\nu}_1 = (m_1,-n_1),\dots,\tilde{\nu}_r = (m_r,-n_r)$. 
We also need to consider the top, left and bottom edges of $P_2$,
which have the following outward-pointing normal vectors:  
\[
\tilde{\nu}_{\text{\it top}} = (0,b),\ \tilde{\nu}_{\text{\it left}} =
(-n,0),\ \tilde{\nu}_{\text{\it bottom}} = (0,-c),
\]
where $n$ is the degree of $\theta(x,y)$ in $y$, $b$ is the degree in
$x$ of the coefficient of $y^n$, and $c$ is the degree of $x$ of the
coefficient of $y^0$.  Note that $n$ and $b$ are the same as in
Minding's paper.

Since $P_1$ lies in the first quadrant and contains the origin, and
$f(x,y)$ has degree $m$ in $y$, we easily obtain
\[
\max_{u\in P_1}\, (u\cdot \tilde{\nu}_{\text{\it top}}) = mb,\
\max_{u\in P_1}\, (u\cdot \tilde{\nu}_{\text{\it left}}) = \max_{u\in
P_1}\, (u\cdot \tilde{\nu}_{\text{\it bottom}}) = 0.
\]
Using \eqref{mfor}, it follows that the mixed area is given by
\begin{equation}
\label{mformin}
\mathcal{M}(P_1,P_2) = mb + \sum_{i=1}^r \max_{u \in P_1} \big(u\cdot
(m_i,-n_i)\big).
\end{equation}
We can interpret the sum in \eqref{mformin} as follows.  Recall that
Minding took the $n$ solutions $y = c_{i,j} x^{-\frac{n_i}{m_i}} +
\cdots,\ i = 1,\dots,r,\ j = 1,\dots,m_i$ from \eqref{gensol} and
substituted them into $f(x,y)$.  If $f(x,y) = \sum_{(k,l) \in
P_1\cap\mathbb{Z}^2} A_{k,l}x^k y^l$, then:
\begin{align*}
\deg\big(f(x,c_{i,j} x^{-\frac{n_i}{m_i}} + \cdots)\big) &=
\deg\Big(\sum_{(k,l) \in P_1\cap\mathbb{Z}^2} A_{k,l}x^k (c_{i,j}
x^{-\frac{n_i}{m_i}} + 
\cdots)^l\Big) \\
&= \deg\Big(\sum_{(k,l) \in P_1\cap\mathbb{Z}^2}
A_{k,l}c_{i,j}x^{k+l(-\frac{n_i}{m_i})} + 
\cdots\Big)\\ 
&= \max_{(k,l) \in
P_1\cap\mathbb{Z}^2}\big(k+l(-\textstyle{\frac{n_i}{m_i}})\big), 
\end{align*} 
where the final equality holds for generic coefficients.  Doing
this for each of the $n$ solutions gives the $n$ numbers that
Minding denoted $k_1,\dots,k_n$.  Then the sum of the first $m_1$
of these numbers is
\begin{align*}
m_1\,\max_{(k,l) \in P_1\cap\mathbb{Z}^2}
\big(k+l(-\textstyle{\frac{n_1}{m_1}})\big) 
&= \max_{(k,l) \in P_1\cap\mathbb{Z}^2}\big(m_1k+l(-n_1)\big)\\ 
&= \max_{(k,l) \in P_1\cap\mathbb{Z}^2}
\big((k,l)\cdot(m_1,-n_1)\big). 
\end{align*}
This is clearly the first term of the sum in \eqref{mformin}.
The other terms match up in exactly the same way, which gives the
formula
\[
\mathcal{M}(P_1,P_2) = mb + k_1 + \dots + k_n.
\]
This proves that Minding's formula is the mixed area of the 
Newton polygons. 

\subsection{Minding's Formula and Mixed Subdivisions} 
\label{sub:mix} 
Here we present an alternative approach to seeing how 
Minding's formula relates to mixed area. The key idea is a 
combinatorial construction due to Huber and Sturmfels \cite{hs}. 

For polygons $P_1,P_2 \subset \R^2$, a \emph{mixed
subdivision}\footnote{Mixed subdivisions are usually defined in terms
of $n$-tuples of polytopes in $\R^n$.  For simplicity, our definition
for $n = 2$ is more restrictive than usual.  A complete exposition on
mixed subdivisions can be found in \cite{rojas} in this volume or
\cite{hs}.}  is a collection of polygons $\cC = \{C_i\}^N_{i=1}$,
called \emph{cells}, satisfying the following conditions:\index{mixed 
subdivision} 

\goodbreak

\begin{coxenum} 
\item $\bigcup_{C_i \in \cC} C_i = P_1+P_2$.
\item For any two distinct $C_i$ and $C_j$ in $\cC$,  $C_i\cap C_j$ is a 
$($possibly empty\/$)$ face of both $C_i$ and $C_j$. 
\item There are two distinguished cells $C_i$ and $C_j$, called 
\emph{unmixed cells}, such that  
$C_i = P_1+p$ and $C_j = P_2+q$ for some $p,q \in \R^2$. 
\item All remaining cells, called the \emph{mixed
cells}, are parallelograms.  
\end{coxenum}

\begin{ex} 
{\rm Taking $P_1$ and $P_2$ to respectively be the Newton polygons 
of Minding's $f$ and $\theta$, one possible mixed subdivision 
is the following:\index{mixed area}\index{mixed subdivision}}

\smallskip

\[
\epsfig{file=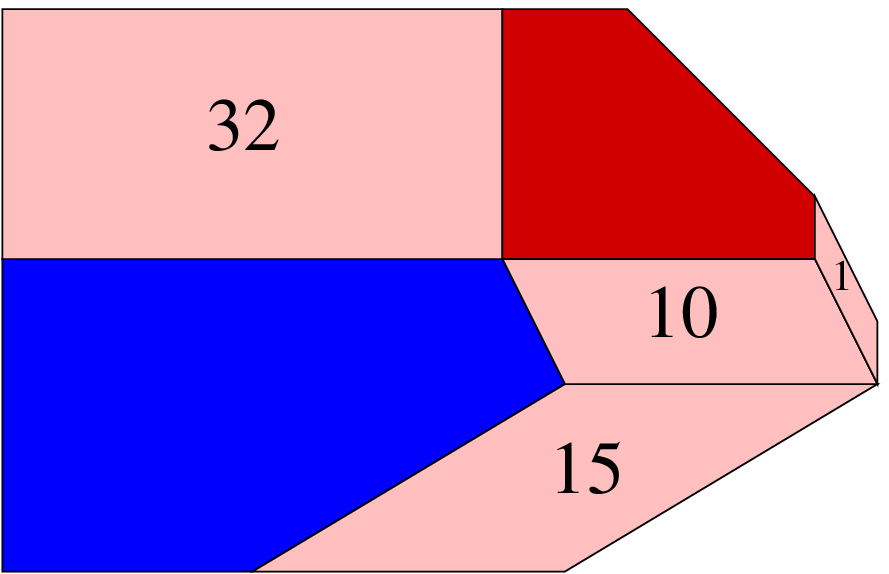,height=1.5in}
\]

\vskip-10pt

\hfill\qed
\end{ex}

Note in particular that a mixed subdivision of the Minkowski sum
$P_1+P_2$ allows us to express the Inclusion-Exclusion formula for
mixed area as a sum of parallelogram areas. For instance, the mixed
area from Minding's example is the sum of the parallelogram areas
above:\ $32+10+1+15 = 58$.

We can identify the individual terms of Minding's formula as polygon
areas as follows.  Let $P_1$ and $P_2$ be lattice polygons of the form
\eqref{genpoly} from Section \ref{sub:puiseux}.  The precise meaning
of this is that $P_1$ and $P_2$ lie in the first quadrant, contain the
origin, and satisfy
\[
\max_{(s,t) \in P_i}(t) = \mathrm{Length}(P_i\cap\{x = 0\}),\quad i =
1,2.
\]
Define
\begin{align*}
m &= \mathrm{Length}(P_1\cap\{x = 0\})\\
n &= \mathrm{Length}(P_2\cap\{x = 0\})\\
b &= \mathrm{Length}(P_2\cap\{y = n\}) 
\end{align*}
and let $F_i$, $i = 1,\dots,r$ be the edges of $P_2$ corresponding to
the vectors $v_i = (n_i,m_i)$ in \eqref{genpoly}.  Using \cite[Sec.\
7]{rojas}, one can show that there is a mixed subdivision $\cC$ of
$P_1$ and $P_2$ such that:

\smallskip

\begin{coxenum}
\item $P_1+(b,n)$ and $P_2 + (0,0)$ are the unmixed cells.
\item $[0,b]\times [n,n+m]$ is a mixed cell of area $mb$.
\item All other mixed cells are of the form $E+F_i$ for certain edges
$E$ of $P_1$ and $i\!\in\!\{1,\ldots,r\}$.
\item For any $i \!\in\!\{ 1,\dots,r\}$, let $\widetilde{C}_i$ be the union of
all mixed cells of $\cC$ of the form $E+F_i$ where $E$ is an 
edge of $P_1$.  Then there is an
rational number $\ell_i$ such that the parallelogram $[0,\ell_i] +
F_i$ has the same area as $\widetilde{C}_i$.  Furthermore,
\[
\mathrm{Area}(\widetilde{C}_i) = \mathrm{Area}([0,\ell_i] + F_i) = m_i
\ell_i = \max_{u \in P_1} \big(u\cdot(m_i,-n_i)\big),
\]
where $\max_{u \in P_1} \big(u\cdot(m_i,-n_i)\big)$ is the $i$th term
of the sum in \eqref{mformin}.
\end{coxenum}

The idea is that the parallelogram $[0,\ell_i] + F_i$ ``straightens
out'' $\widetilde{C}_i$.  For instance, in the subdivision of Example
1, let $F_1$ be the edge of $P_2$ corresponding to $v_1 = (-1,2)$.
Since $F_1$ has height $2$, the two mixed cells in Example 1 involving this
edge give the following polygons of the same area:
\[
\begin{picture}(250,35)
\thicklines
\put(5,25){\line(1,0){50}}
\put(15,5){\line(-1,2){10}}
\put(15,5){\line(1,0){50}}
\put(65,5){\line(-1,2){10}}
\put(65,5){\line(0,1){10}}
\put(55,25){\line(0,1){10}}
\put(65,15){\line(-1,2){10}}
\put(80,12){$\longrightarrow$}
\put(105,25){\line(1,0){50}}
\put(115,5){\line(-1,2){10}}
\put(115,5){\line(1,0){50}}
\put(165,5){\line(-1,2){10}}
\put(165,5){\line(1,0){5}}
\put(155,25){\line(1,0){5}}
\put(170,5){\line(-1,2){10}}
\put(-3,10){{\small $F_1$}}
\put(40,-3){{\small 5}}
\put(68,7){{\small 1}}
\put(140,-3){{\small 5}}
\put(165,-4){{\small $\frac12$}}
\put(200,10){$\ell_1 = 5+\frac12 = \frac{11}2$}
\put(35,12){10}
\put(59.6,16){{\tiny 1}}
\put(135,12){10}
\put(161.5,12){{\tiny 1}}
\end{picture}
\] 

\medskip

\noindent In general, comparing part (d) to the formulas at the end
of Section~\ref{recursion}, we see that $\ell_i$ is the degree of
$f(x,c_{i,j}x^{-\frac{n_i}{m_i}} + \cdots)$ for $j = 1,\dots,m_i$.  It
follows that
\[
\qquad\quad
\underbrace{\ell_1,\dots,\ell_1}_{m_1\,
  \text{times}},\ \dots\ ,  
\underbrace{\ell_r,\dots,\ell_r}_{m_r\,
  \text{times}}
\]
are Minding's numbers $k_1,\dots,k_n$.  Thus (a)--(d) imply that
\begin{align*}
\mathcal{M}(P_1,P_2) &= \text{sum of areas of mixed cells in the mixed
  subdivision}\\
&= mb + \sum_{i=1}^r m_i \ell_i = mb + k_1 + \dots + k_n.
\end{align*}
so that Minding's formula indeed gives the mixed area.  

\begin{ex}  
{\rm Here is an example more complicated than Minding's:
\[
\epsfig{file=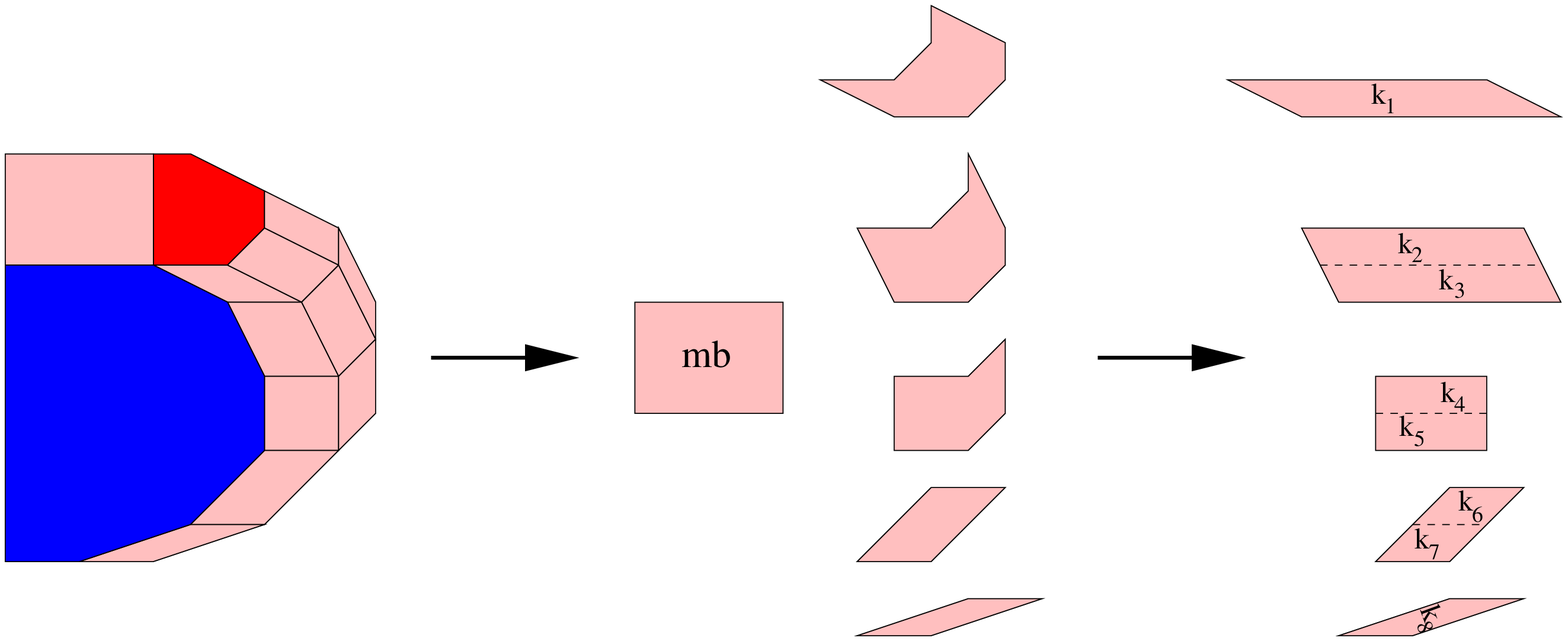,height=1.8in}
\]
Note that $k_2 = k_3$, $k_4 = k_5$, and $k_5 = k_6$.  As noted above,
the right-hand side of our illustration above shows how we can
calculate the area of an entire ``strip'' of mixed cells as the area
of a single parallelogram.} \qed
\end{ex} 

In general, the connections between algebraic geometry and combinatorics 
are very deep. In our context, this can be made into a ``mantra'' of 
sorts: any algorithm for computing mixed area yields a proof of 
Bernstein's Theorem in the plane, and vice-versa. 

\subsection{Affine Roots and Genericity} 
\label{sub:gen} 
Minding's formula reduces to the BKK bound only when we make certain
assumptions on the coefficients.  On the other hand, Minding's formula
is without qualification a special case of the following result due to
Li and Wang:

\begin{thm} \cite{liwang} 
The number of isolated roots in $\C^n$ of a system of 
polynomial equations $f_1(x_1,\ldots,x_n) = \cdots = f_n(x_1,\ldots,x_n)=0$ 
is no more than $\cM(P'_1,\ldots,P'_n)$ where, for all $i$, 
$P'_i$ is the convex hull of the union of the origin and
the Newton polytope of $f_i$. \qed
\end{thm} 

Note in particular that one need not assume that the polynomials have
a nonzero constant term or that their exponents satisfy any
restriction.  That Li and Wang's result holds for generic coefficients 
(assuming that all the polynomials have nonzero constant term) 
follows easily from the BKK bound, so the theorem is interesting
because it states an upper bound which is true regardless of
degeneracies or the presence of constant terms.

However, this brings to mind an important question:\ How much
genericity is necessary for the BKK bound to hold?  From Bernstein's
seminal paper \cite{bernie}, it is easy to show that it suffices to have all
coefficients corresponding to Newton polytope vertices be
generic.\footnote{It is certainly not enough to have all the vertex
coefficients simply be nonzero.}  That one could get away with less
genericity appears to have been first observed in \cite{cr}, and 
a complete characterization of which sets of coefficients need to 
be generic (for the BKK bound, the Li-Wang bound, and even more 
general bounds) appears in \cite{toric}.
Interestingly, although he refrains from explicitly stating a clear
characterization of which coefficients should be chosen generically,
Minding did appear to understand that not all coefficients need to be
generic for his formula to be true.

\subsection{Final Comments} 

For a nice discussion of how Minding's argument relates to the BKK
bound, we refer the reader to Khovanskii's essay in \cite[Ch.\ 4,
Sec.\ 27, Addendum 3]{bz}.  Also, \cite[Sec.\ 6]{rojas} explains
how mixed volume arisely naturally when considering the solutions of
polynomial equations.  

\section*{Acknowledgements} 

We would like to thank Hal Schenck for his comments on the
translation and Ms.\ Nishanti for assistance in generating  
a LaTeX version of the original German version of \cite{minding}. The 
second commentator also acknowledges the support of the 
Hong Kong University Grants Council through grant \#9040402-730.

\end{document}